\documentclass[twoside,12pt]{article}
\usepackage{amsmath}
\usepackage{amssymb}
\usepackage{amscd}
\usepackage{amsthm}%+

\numberwithin{equation}{subsection}
 
\theoremstyle{plain}

\theoremstyle{remark}
\newtheorem{remark}[equation]{Remark}
\theoremstyle{definition}

\newcommand{\R}{\mathbb R}
\newcommand{\C}{\mathbb C}
\def\ga{\gamma}

\def\ra{\rightarrow}

\def\e{\emph}
\def\i{\infty}
\def\p{\partial}
\def\b{\begin}
\newcommand{\ol}{\overline}

\begin{document}

\title{    \flushleft{\bf{Quasim\"obius maps preserve uniform domains}}      }
%\author{\flushleft{Xiangdong Xie}}
\date{  }
\maketitle

\vspace{-5mm}

\noindent 
Xiangdong Xie\newline
Department of Mathematical Sciences,
 University  of Cincinnati,\newline
Cincinnati, OH 45221, U.S.A.\newline
Email:   {  xxie@math.uc.edu}

%\noindent
%{\bf{Mathematics Subject Classification(2000).}} Primary 53C23, 51F99, 57M20; 
% Secondary  53C70, 57M60.   
% \newline
%{\bf{Keywords.}} Tits boundary, Tits metric, CAT(0), 2-complex, geodesic.

%\keywords    test  \endkeywords
%\subjclass  51  \endsubjclass

\pagestyle{myheadings}

\markboth{{\upshape Xiangdong Xie}}{{\upshape  Quasim\"obius maps preserve uniform domains}}

%\footnotetext{test}
      
%\keywords    test  \endkeywords
%\subjclass  51  \endsubjclass
\vspace{3mm}

\noindent
{\small {\bf Abstract.}
We show  that  if  a domain  $\Omega$ in a 
geodesic metric space is 
  quasim\"obius   to a uniform domain in some metric space, then $\Omega$
 is also uniform.}

%Your abstract number is: 1020-30-15.

%\tableofcontents  

\vspace{3mm}
\noindent
{\small {\bf{Mathematics Subject Classification (2000).}}   30C65
%20F65. %57M20, 20F67, 20E07.
}

%30C65  (1991-now) Quasiconformal mappings in ${\bf R}^n$, other generalizations
%20F34 Fundamental groups and their automorphisms 
%20F65 Geometric group theory 
%20F67 Hyperbolic groups and nonpositively curved groups
%20F99 None of the above, but in this section
%20E07 Subgroup theorems; subgroup growth
%53C20 Global Riemannian geometry, including pinching
%53C23 Global topological methods 
%57M20 Two-dimensional complexes
%57M60 Group actions in low dimensions
%(1) (2) (a) (b) (i) (ii)
%20F69 Asymptotic properties of groups

%$\{x_i\}_{i=1}^\i$  converges to $\xi\in \ol{X}$ 
% an index two subgroup

\vspace{3mm}

%hyperbolic element, parabolic element, 
%quasi-convex.}

\setcounter{section}{1}
\setcounter{subsection}{0}

\subsection{Introduction}\label{s1}

In this paper we study the behavior of uniform domains under quasisymmetric and 
quasim\"obius maps between metric spaces.

Let $(X, d)$ be 
a metric space.   
A subset of $X$ is called a \e{domain} if it is open and connected.
 Let  $0<\lambda\le 1$ and $c\ge 1$.
   We say a domain $\Omega\subset X$ 
with $\p\Omega\not=\emptyset$
is \e{$(\lambda, c)$-quasiconvex}, if for any $x\in \Omega$, 
  and any two points $y_1, y_2\in B(x, \lambda d(x, \p\Omega))$, there is a 
path $\gamma$ in $\Omega$ 
from $y_1$ to $y_2$ with  length $\ell(\gamma)\le c\, d(y_1, y_2)$.

We say $(X, d)$  is \e{$c$-quasiconvex}
if  for any two points $x, y\in X$, there is a path $\gamma$ joining $x$ and $y$ 
 with length $\ell(\gamma)\le c\, d(x,y)$.    A  metric space $(X, d)$ is 
quasiconvex  if it is $c$-quasiconvex for some $c\ge 1$.
 Notice that if $(X,d)$ is $c$-quasiconvex, then 
 every domain  in $X$  is $(\lambda, c)$-quasiconvex for 
   all  
$0<\lambda\le 1$. 

Recall that   a  metric space is \e{proper} if all its closed balls  are compact.
%We prove the following theorem in this paper.
 The main result of the paper is as follows:

\b{Th}\label{main}
{Let $(X_i, d_i)$ \e{(}$i=1, 2$\e{)} be a proper metric space,  and 
     $\Omega_i\subset X_i$ a domain with $\p\Omega_i\not=\emptyset$.  
  Suppose $(\Omega_1, d_1)$ is $c_1$-uniform,    
$\Omega_2\subset (X_2, d_2)$ is $(\lambda, c_2)$-quasiconvex
 for some $0<\lambda\le 1$ and $c_2\ge 1$,   and   
    there is a  $\eta$-quasim\"obius homeomorphism
$ (\Omega_1, d_1)\ra (\Omega_2, d_2)$.   %\newline 
    Then 
$(\Omega_2, d_2)$ is  
  $c$-uniform  for some constant $c$.  %\newline
%(2)  If  $(\Omega_2, d_2)$  is unbounded,  then $(\Omega_2, d_2)$ is  
   %$c'$-uniform, where $c'$ depends only on  $\lambda$,  $c$  and $c_1$.
}

\end{Th}

%Let  $c\ge 1$. %  be  a constant. 

It is reasonable to  assume  %some kind of 
quasiconvexity in Theorem \ref{main}, at least 
one needs that the domains are rectifiably connected, as shown by the following example
(The author thanks  David  Herron
%and Nageswari Shanmugalingam
 for pointing out this example).
Let $\Delta\subset \R^2$ be the open unit disk in the plane with the Euclidean metric
 $d$. 
The domain $\Delta$ is  clearly uniform.  For any $0<\epsilon<1$, the identity map 
$(\Delta, d)\ra (\Delta, d^\epsilon)$ is $\eta$-quasim\"obius with
$\eta(t)=t^\epsilon$. But $(\Delta, d^\epsilon)$ is not uniform since there are    no rectifiable curves in $(\Delta, d^\epsilon)$  except the constant curves.

Theorem \ref{main}   is not  quantitative in the  sense   that  
   when $(\Omega_2, d_2)$ is bounded 
one can not control the constant $c$ in terms of $c_1$, $c_2$, $\lambda$ 
  and  $\eta$  alone.  See Section \ref{s5}  for an example. 
 It is unclear whether one can make Theorem \ref{main}  
quantitative  under the  stronger  assumption that 
$(\Omega_2, d_2)$  is   quasiconvex.  See Section \ref{s5}
    for more detail about this question.
%Theorem 
%\ref{main} says that if 
%$(\Omega_1, d_1)$ is 
%$c_1$-uniform,   $\Omega_2\subset (X_2, d_2)$ is $(\lambda, c_2)$-quasiconvex
 %and $(X_2, d_2)$ is $c_2$-quasiconvex, 
%and 
%there is a $\eta$-quasim\"obius   map
%$g: \Omega_1\ra \Omega_2$, then 
%$(\Omega_2, d_2)$ is  
%$c$-uniform for some constant $c\ge 1$.
%However, it is unclear whether 
%one  can  control $c$ in terms of $\lambda$,  $c_1$, $c_2$ and 
%$\eta$  alone.
%See the example in Section \ref{s4}. 
 On the other hand, when $(\Omega_2,  d_2)$ is unbounded  or 
  when $(X_2,  d_2)$  is    quasiconvex and annular convex, 
 there is the following   quantitative result. Recall that a metric space is 
\e{$c$-annular convex} for some $c\ge 2$
if for any  
   $x\in X$,  each  $r>0$, and any 
   $y, z\in B(x,   2r)\backslash B(x, r)$, 
  there is   a path $\gamma$ from $y$ to $z$  satisfying  
 $\ell(\gamma)\le c \, d(y, z)$   and 
    $\gamma \cap B(x, r/c)=\emptyset$.

\b{Th}\label{ac}  %\e{(\cite{HSX})}
{Let $(X_i, d_i)$ \e{(}$i=1, 2$\e{)} be a proper metric space,  and 
     $\Omega_i\subset X_i$ a domain with $\p\Omega_i\not=\emptyset$.
Suppose   $(\Omega_1, d_1)$ is $c_1$-uniform, and  %$(X_2, d_2)$ is $c_2$-quasiconvex and  
  %$c_2$-annular convex,  and 
there is an   $\eta$-quasim\"obius homeomorphism
$h:   (\Omega_1, d_1)\ra (\Omega_2, d_2)$.     \newline
 (1)    If    $(\Omega_2, d_2)$  is unbounded  and
$(\lambda, c_2)$-quasiconvex,  then $(\Omega_2, d_2)$ is
$c$-uniform with $c=c(\eta, c_1, c_2, \lambda)$.
   \newline  % then $(\Omega_2, d_2)$ is  
(2)  If    $(X_2, d_2)$ is $c_2$-quasiconvex and  
$c_2$-annular convex,
   %$c'$-uniform, where $c'$ depends only on  $\lambda$,  $c$  and $c_1$.
%and a domain $\Omega$ 
%in $X$  with $\p\Omega\not=\emptyset$ 
%is $\eta$-quasim\"obius to a  $c_2$-uniform domain in some  metric space. 
  then $(\Omega_2, d_2)$ is $c$-uniform with $c=c(\eta, c_1, c_2)$.}

\end{Th}

Theorem \ref{ac} (2)  was  first  proved in \cite{HSX}
by using a characterization of uniform domains in terms of Gromov hyperbolic spaces and 
 the quasiconformal structure on the Gromov boundary.  In this paper we give a different proof.

For quasisymmetric maps, we have the following:

\b{Th}\label{t4.1}
%{\bf{Theorem \ref{t4.1}}}.
{Let $(X_i, d_i)$ \e{($i=1, 2$)} be a proper metric space,   
     $\Omega_i\subset X_i$ a domain
 with $\p\Omega_i\not=\emptyset$,  and $g: (\Omega_1, d_1) \ra (\Omega_2, d_2)$
  an    $\eta$-quasisymmetric map.  Suppose $(\Omega_1, d_1)$ is 
$c_1$-uniform  and $(\Omega_2, d_2)$ is $ (\lambda, c_2)$-quasiconvex.  Then 
$(\Omega_2, d_2)$ is 
$c$-uniform  with $c=c(\eta,   c_1, c_2, \lambda)$. }

\end{Th}

In the case when the metric spaces are geodesic spaces 
Theorem \ref{t4.1} essentially has been proved by V\"ais\"al\"a
( see the proofs of Lemma 10.21 and Theorem 10.22 in  \cite{V}), who stated it only for 
 Banach spaces. 
%is a slight improvement of an result of V\"ais\"al\"a (??\cite{V}), 
%%who 
 %essentially  proved Theorem \ref{t4.1}  in the case when $(X_2, d_2)$ is quasiconvex.
 Our proof is  a slight modification of V\"ais\"al\"a's. 
%The reason that 
We include  this theorem  here since  we need it for the proofs of
 Theorem \ref{main} and Theorem \ref{ac}. 

The invariance of uniform domains under quasim\"obius  maps 
 was  (implicitly) obtained by 
  Gehring  and    Martio (\cite{GM}) for Euclidean domains, and has been 
  established by 
V\"ais\"al\"a  for domains in Banach spaces (\cite{V2}, \cite{V}).

Theorem \ref{main} is proved by using Theorem \ref{t4.1} and a construction of 
Bonk-Kleiner.
For any unbounded proper  metric space $(X, d)$ and $p\in X$, 
   Bonk and  Kleiner constructed  a metric $\hat{d}_p$ on the 
 one point compactification $X\cup \{\i\}$ of $X$, such that the identity map 
$(X, d)\ra (X, \hat{d}_p)$  is quasim\"obius (see \cite{BK} or Section \ref{s2}).  
Furthermore, for any domain
$\Omega\subset X$ with $\p\Omega\not=\emptyset$,   % with $p\in \partial \Omega$,  
$\Omega$ is uniform 
 with respect to $d$ if and only if 
$\Omega$ is  uniform 
 with respect to $\hat{d}_p$ (see Theorem \ref{bk2}).

Since a quasim\"obius map between bounded metric spaces is quasisymmetric, 
 Theorem \ref{main} follows from Theorem \ref{t4.1} when both
 $(\Omega_1, d_1)$ and  $(\Omega_2, d_2) $ are bounded.
 In general, when  we consider the metric  $d_i'=:\hat{d_i}_{p_i}$ ($p_i\in X_i$)
   on $X_i$ ($i=1, 2$),
  a quasim\"obius map 
$g: (\Omega_1, d_1) \ra (\Omega_2, d_2)$  becomes a quasisymmetric 
map $g: (\Omega_1, d_1') \ra (\Omega_2, d_2')$ as
 $(X_i, d_i')$
 are   bounded  for  $i=1,  2$.  Theorem \ref{t4.1} implies that 
   $(\Omega_2, d_2')$ is uniform, and hence 
$(\Omega_2, d_2)$ is also uniform by the preceding paragraph.

\noindent
{\bf{Acknowledgment}}.
The author would like to thank 
Nageswari Shanmugalingam
 for carefully reading an earlier version and suggesting many improvements.

\subsection{Preliminaries}\label{s2}

In this Section we recall some basic definitions and facts, see \cite{V}
 and \cite{BHX} 
 for more details.

%Here we recall the definition of quasihyperbolic metrics.

Let $(X, d)$ be a complete 
   metric space, and 
 $\Omega\subset X$  a domain.
The metric boundary of $\Omega$ is $\p \Omega=\overline \Omega\backslash \Omega$.
In this paper we always assume 
  $\p\Omega\not=\emptyset$. 
For $x\in \Omega$, we denote $d(x)=d(x, \p \Omega)$.
We say $\Omega$ is \e{rectifiably connected}
 if for any $x,y\in \Omega$ 
  there is a path in $\Omega$ from $x$ to $y$ with finite length. 
For a rectifiably connected domain $\Omega$, 
 %The metric boundary of $\Omega$ is $\p \Omega=\overline \Omega-\Omega$.
%For $x\in \Omega$, we denote $d(x)=d(x, \p \Omega)$.  
  the \e{quasihyperbolic metric }
$k$ on $\Omega$ is defined as follows:
for $x, y\in \Omega$, 
$$k(x,y):=\inf \int_\gamma\frac{1}{d(z)}\,ds(z),$$  
where  %$|dz|$ denotes the arc length element along $\ga$ and 
$\gamma$  runs over all rectifiable  curves in $\Omega$ joining $x$ and $y$.
 Here $ds$ denotes the arc length element along $\ga$.
For $x,y\in \Omega$, we set 
$$r_\Omega(x,y)=\frac{d(x,y)}{d(x)\wedge d(y)} \;\;\;\;{\text{and}}\;\;\;
j_\Omega(x,y)=\log\left(1+r_\Omega(x,y)\right),$$ 
where $a\wedge b$ denotes $\min\{a,b\}$ for real numbers $a, b$.

The length metric on $\Omega$ is defined  as follows: for $x,y\in \Omega$, 
$l_\Omega(x,y)$ is the infimum of length of paths  in $\Omega$ from 
 $x$ to $y$.

\b{Prop}\label{homeosamet} \e{(Proposition 2.8 of \cite{BHK})}
{Suppose $(\Omega, d)$ is locally compact and rectifiably connected.  
If the identity map ${\text{id}}:  (\Omega, d)\ra (\Omega, l_\Omega)$ is 
 a homeomorphism, then ${\text{id}}:  (\Omega, d)\ra (\Omega, k)$
  is also a homeomorphism, %  $(\Omega, d)\ra (\Omega, k)$, 
  and $(\Omega, k)$  is a  proper  geodesic space.}

\end{Prop}

We observe that if $(\Omega, d)$ is $(\lambda, c)$-quasiconvex for some
 $0<\lambda\le 1$ and $c\ge 1$,  then 
the identity map $(\Omega, d)\ra (\Omega, l_\Omega)$ is 
locally bilipschitz and hence is 
 a homeomorphism.   % Notice that 
   However,   
${\text{id}}:   (\Omega, d)\ra (\Omega, l_\Omega)$   is not always  a 
homeomorphism ( one can easily construct an example using topologist's sine curve).

\b{Le}\label{bal} \e{(Theorem 3.7 (1)  of \cite{V})}
 {The following holds for  all  $x, y\in \Omega$,  
\[ k(x,y)\ge j_\Omega(x,y)\ge \log\frac{d(x)}{d(y)}.\]
 % for all
  % $x,y\in \Omega$}
}
\end{Le}

%Let $(X,d)$ be a metric space and $\Omega\subset X$ an open subset
 % with $\p\Omega\not=\emptyset$.
  Let $c\ge 1$.    A path $\gamma:[0,1]\ra \Omega$
is called a \e{$c$-uniform curve}
if:\newline
(1) $\ell(\gamma)\le c\;d(\gamma(0), \gamma(1))$;\newline
(2) $c\;d(\ga(t))\ge \ell(\gamma|[0,t])\wedge \ell(\gamma|[t,1])$ for all $t\in [0,1]$.\newline
The domain 
$\Omega\subset (X,d)$ is called a \e{$c$-uniform domain} 
  in $(X,d)$  
if  every  two points $x, y\in \Omega$  
can be joined by a
 $c$-uniform curve.

Let $\eta: [0,\i)\ra [0,\i)$ be a homeomorphism.
 A homeomorphism between metric spaces
$f: (X, d_X)\ra (Y, d_Y)$ is \e{$\eta$-quasisymmetric} if for all pairwise 
 distinct points $x, y, z\in X$, we have
$$\frac{d_Y(f(x), f(y))}{d_Y(f(x), f(z))}\le \eta\left(\frac{d_X(x,y)}{d_X(x,z)}\right).$$
A homeomorphism 
$f: (X, d_X)\ra (Y, d_Y)$ is
quasisymmetric  if it is $\eta$-quasisymmetric  for some $\eta$.
The inverse of a quasisymmetric map is quasisymmetric, and the composition of two 
quasisymmetric maps is also quasisymmetric.

Let  $Q=(x_1, x_2, x_3, x_4)$ be  a quadruple of pairwise
distinct points in $(X, d)$.
The \e{cross ratio} of $Q$ with respect to the metric $d$    is:
$$cr(Q, d)=\frac{d(x_1, x_3)d(x_2,  x_4)}{d(x_1,  x_4) d(x_2,  x_3)}.$$
Let $\eta: [0, \infty)\ra [0, \infty)$ be a homeomorphism. 
A homeomorphism between metric spaces $f: (X, d_X)\ra (Y,  d_Y)$ is 
  an    \e{$\eta$-quasim\"obius} map if 
$$cr(f(Q), d_Y)\le \eta(cr(Q, d_X))$$
for all
 quadruple $Q$ of distinct points in $X$, where $f(Q)=(f(x_1), f(x_2), f(x_3), f(x_4))$. 
  A homeomorphism 
$f: (X, d_X)\ra (Y, d_Y)$ is  quasim\"obius
  if it is $\eta$-quasim\"obius  for some $\eta$.
The inverse of a quasim\"obius map is quasim\"obius, and the composition of two quasim\"obius
 maps is also quasim\"obius.

Quasisymmetric  maps are quasim\"obius,
 and quasim\"obius  maps between bounded metric spaces are quasisymmetric.

Let $(X,d)$  be an unbounded metric space and $p\in X$.  Set $S_p(X)=X\cup \{\i\}$,
 where $\i$ is a point not in $X$.
 %  be the one point compactification of $Z$. 
    Let  
\[s_p(x,y)=\frac{d(x,y)}{(1+d(x,p))(1+d(y,p))}\]
for  $x, y\in X$,   $s_p(x,\i)=s_p(\i,x)=\frac{1}{1+d(x,p)}$  
 for $x\in X$  and $s_p(\i,\i)=0$.
%Let $\hat{h}_p:  S_p(X)\ra [0,\i)$ be defined as follows:
%$\hat{h}_p(x)=\frac{1}{1+d(x,p)}$ if $x\in X$  and $\hat{h}_p(\i)=0$.  
%For $x,y\in X$, let $s_p(x,y)=\hat{h}_p(x)\hat{h}_p(y)d(x,y)$, 
 %  $s_p(x,\i)=s_p(\i,x)=\hat{h}_p(x)$  and $s_p(\i,\i)=0$.  
 For $x, y\in S_p(X)$,
   define 
$$\hat{d_p}(x, y):=\inf \Sigma_{i=0}^{k-1}s_p(x_i, x_{i+1}),$$
where the infimum is taken  over all finite sequences of points  
$x_0, \cdots, x_k\in S_p(X)$
   with  $x_0=x$  and $x_k=y$.  Then  $\hat{d_p}$ is a metric on $S_p(X)$  
 and 
$$\frac{1}{4}s_p(x,y)\le \hat{d_p}(x,y)\le s_p(x,y)\;\;\; 
{\text{for}}\;\;\;  x,y\in S_p(X).$$
%whose  induced topology 
 %on $X-\{a\}$  
%agrees with 
 % the  subspace topology  on $X-\{a\}\subset X$. 
   Furthermore   
 the identity map 
${\text{id}}:  (X, d)\ra (X, \hat{d_p})$   is an $\eta$-quasim\"obius  homeomorphism  
 with   $\eta(t)=16t$.
If  $(X, d)$
is  $c$-quasiconvex and $c$-annular convex, then 
$(S_p(X), \hat{d}_p)$  is  $c'$-quasiconvex and  $c'$-annular convex,
  where $c'$ depends only on $c$.    See \cite{BHX}  for a proof of the above 
 statements.

\b{Th}\label{bk2}\e{(\cite{BHX})}
{Let $(X, d)$ be an unbounded proper   metric space,   $\Omega\subset X$ 
a  %unbounded 
%rectifiably connected
 domain
with $\p \Omega\not=\emptyset$   
  and $p\in X$.  Then 
$\Omega\subset (X, d)$ is  uniform  if and only if 
$\Omega\subset (S_p(X), \hat{d}_p)$ is uniform. 
Furthermore,\newline
 (1)  if $(\Omega, d) $  %\subset X$ 
  is unbounded, $p\in \p\Omega$   and $\Omega\subset (X, d)$ is  $c$-uniform,
 then $\Omega\subset (S_p(X), \hat{d}_p)$ is $c'$-uniform
with $c'$ depending only on $c$;\newline
(2) if $(\Omega, d) $  %\subset X$ 
  is unbounded  and $\Omega\subset (S_p(X), \hat{d}_p)$ is $c$-uniform, then 
$\Omega\subset (X, d)$ is  $c'$-uniform with $c'$ depending only on $c$.
}

\end{Th}

Let $(X, d)$ be a metric space and $p\in X$.  
Set $I_p(X)=X\backslash\{p\}$ if $X$ is bounded and   $I_p(X)=(X\backslash\{p\})\cup \{\i\}$ 
if $X$ is unbounded, 
   where  $\i$ is a point not in  $X$. 
We shall  define a metric $d_p$ 
on  $I_p(X)$.

%Define  a function $h_p:  I_p(X)\ra [0,\i)$ by 
 %$h_p(x)=\frac{1}{d(x,p)}$  for   $x\in X\backslash\{p\}$  and $h_p(\i)=0$.   

 Let 
\[f_p(x,y)=\frac{d(x,y)}{d(x,p)d(y,p)}\]   
for $x,y\in X\backslash\{p\}$,  
  $f_p(x,\i)=f_p(\i, x) =\frac{1}{d(x,p)}$ for $x\in X\backslash\{p\}$
  and  $f_p(\i,\i)=0$.

For $x, y\in I_p(X)$,
  we define 
$$d_p(x, y):=\inf \Sigma_{i=0}^{k-1}f_p(x_i, x_{i+1}),$$
where the infimum is taken  over all finite sequences of points  
$x_0, \cdots, x_k\in I_p(X)$
   with  $x_0=x$  and $x_k=y$.
Then
the following holds for all $x,y\in I_p(X)$: 
$$\frac{1}{4}f_p(x,y)\le d_p(x,y)\le f_p(x,y). \leqno{(2.1)}$$
 In particular,  
  $d_p$ is a metric on $I_p(X)$  
%whose  induced topology 
 %on $X-\{a\}$  
%agrees with 
 % the  subspace topology  on $X-\{a\}\subset X$. 
  and  
 the identity map 
${\text{id}}:  (X\backslash\{p\}, d)\ra (X\backslash\{p\}, d_p)$   is an $\eta$-quasim\"obius  homeomorphism  
 with   $\eta(t)=16t$.
If  $(X, d)$
is  $c$-quasiconvex and $c$-annular convex, then the space 
$(I_p(X), d_p)$  is  $c'$-quasiconvex and $c'$-annular convex, where $c'$ depends only on $c$.  
See \cite{BHX}  for a proof of the above 
 statements.

% (see \cite{BHX}).

\b{Th}\label{bk3}\e{(\cite{BHX})}
{Let $(X, d)$ be a  proper  metric space,    $\Omega\subset X$  a  domain 
 %rectifiably connected
%open subset  %with $\p \Omega\not=\emptyset$
  and $p\in \p\Omega$. Assume $\p\Omega$ contains at least two points if $(\Omega, d)$ is bounded. Then 
$\Omega\subset (X, d)$ is  uniform  if and only if 
$\Omega\subset (I_p(X), {d}_p)$ is uniform. Furthermore,  
\newline
(1) if $\Omega\subset (X, d)$ is  $c$-uniform, then 
$\Omega\subset (I_p(X), d_p)$ is $c'$-uniform, where $c'$ depends only on $c$;\newline
(2)  if $X$ is $c$-quasiconvex and $c$-annular convex,
 and $\Omega \subset (I_p(X), d_p)$ 
  is $c_1$-uniform,  then 
$\Omega \subset (X, d)$ is $c'$-uniform with $c'=c'(c,c_1)$.}

\end{Th}

The above two constructions, $S_p(X)$ and $I_p(X)$, are in a  sense
inverse to each other,
% in a certain sense, 
as shown by the following result.

\b{Le}\label{rebm1}\e{(\cite{BHX})}
{Let  $(X, d)$ be an unbounded metric space and $p\in X$. Set $Y=S_p(X)=X\cup \{\i\}$ and 
 denote by $d'$ the metric $(\hat{d_p})_\i$ on  $I_\i(Y)=X$. 
%Let $X=Z\cup \{\i\}$ be  the one point compactification with the metric $\hat d_p$,
  %and $d_\i$ the metric on $Z$. 
    Then the identity map 
${\text{id}}: (X, d)\ra (I_\i(Y), d')=(X, d')$ is   16-bilipschitz.}

\end{Le}

\subsection{Equivalent definitions of uniform domains}\label{s3}

In this Section we prove  the equivalence of several 
different definitions of uniform 
domains. We need this in Section \ref{s4}.
The results in this Section  can  be proved by slightly modifying 
 either the proofs
of   Gehring-Osgood (Section 2 of \cite{GO})  or those  of  
 V\"ais\"al\"a (Section 10 of \cite{V}). 
Both V\"ais\"al\"a's  and Gehring-Osgood's
 proofs use the fact that the metric spaces are geodesic, and this 
    is the only place 
 that needs to be modified. 
  We adopt V\"ais\"al\"a's proofs.  
We include them  here mainly for
 completeness.

Let $(X, d)$  be a proper metric space, and 
     $\Omega\subset X$ a rectifiably connected 
domain with $\p\Omega\not=\emptyset$. Recall that we always have 
$k(x,y)\ge j_\Omega(x,y)$ for all  $x,y\in \Omega$. 
Let $c\ge 1$. We say 
$\Omega\subset X$ is  a \e{QH $c$-uniform  domain}  
if  $k(x,y)\le c\; j_\Omega(x,y)$ for all $x, y\in \Omega$.

We first recall several lemmas from \cite{V}.

\b{Le}\label{l2} (Lemma 10.7 in \cite{V})
{Let $\Omega\subset X$ be a QH $c$-uniform  domain,  $r>0$,
  and  $\gamma$ be 
 a quasihyperbolic geodesic in $\Omega$  
 such that $d(z)\le r$ for all $z\in \gamma$.  
  %contained in $N_r(\p \Omega)$.
  Then $\ell(\gamma)\le M_1(c)\;r$, where $M_1(c)$ is a constant depending only on $c$.}

\end{Le}

\b{Le}\label{l3}  (Lemma 10.8 in \cite{V})
{For each $c\ge 1$, there is a number $q=q(c)\in (0,1)$ with the following property:
  Let $\Omega\subset X$  be a QH  $c$-uniform domain, 
   $\gamma\subset \Omega$  a quasihyperbolic
geodesic 
 with endpoints $a_0$, $a_1$, and let $x\in \gamma$  be a point  with $d(x)\le q \;d(a_0)$.
 Then  for $\gamma_x=\gamma[x,a_1]$  we have $\ell(\gamma_x)\le M_2(c)\,d(x)$,
where 
 $M_2(c)$ is a constant depending only on $c$.}

\end{Le}

\b{Le}\label{es}
{Let %$(X, d)$  be a proper metric space, and 
     $\Omega\subset X$ be  a $(\lambda_0, c_0)$-quasiconvex domain  for some
$0<\lambda_0\le 1$  and  $c_0\ge 1$.
Let $x,y\in \Omega$. 
   If  $\frac{d(x,y)}{d(x)}\le \frac{\lambda_0}{2c_0}$, then 
$k(x,y)  \le 2c_0  \frac{d(x,y)}{d(x)}       %r_\Omega(x,y)
\le \lambda_0$.  In particular,  
  if $r_\Omega(x,y)\le \frac{\lambda_0}{2c_0}$, then 
$k(x,y)  \le 2c_0r_\Omega(x,y)
\le \lambda_0$.}

\end{Le}

\b{proof}
%We may assume $d(x)\le d(y)$. 
The assumption implies 
$d(x,y)\le \frac{\lambda_0}{2c_0}d(x)$. Since 
$\Omega\subset X$ is  $(\lambda_0, c_0)$-quasiconvex,
 there is a path $\gamma$ from $x$ to $y$ with 
$\ell(\gamma)\le c_0 d(x,y)\le \frac{\lambda_0}{2}d(x)\le d(x)/2$.
  It follows that $d(z)\ge d(x)/2$ for all $z\in \gamma$. 
 Now 
\[k(x,y)\le \int_\gamma \frac{1}{d(z)}\,ds(z)\le \frac{2}{d(x)}\ell(\gamma)
\le \frac{2}{d(x)}c_0\,d(x,y)  %\le 2c_0    \frac{d(x,y)}{d(x)}    %r_\Omega(x,y)
\le \lambda_0.\]

%\frac{\lambda_0}{2}d(x)=\lambda_0$.

\end{proof}

\b{Th}\label{cigart}
{Let $(X, d)$  be a proper metric space, and 
     $\Omega\subset X$   be 
   a $(\lambda_0, c_0)$-quasiconvex domain  for some
$0<\lambda_0\le 1$  and  $c_0\ge 1$.
If 
 $\Omega\subset X$ is  QH $c_1$-uniform,  % and that $\gamma\subset \Omega$ is a  geodesic with endpoints $a_0$, $a_1$.  
 then $\Omega\subset X$ 
is  $c_2$-uniform  with $c_2=c_2(\lambda_0, c_0, c_1)$.
%Then there is $c_1=c_1(c)$ such that \newline 
%(1);;;\newline
%(2) $l(\gamma)\le c_1 d(a_0, a_1).$
}
\end{Th}

\b{proof} Fix $a_0, a_1\in \Omega$ and let $\gamma$ be a 
quasihyperbolic geodesic from $a_0$ to  $a_1$.  We shall show that 
$\gamma$ is a $c_2$-uniform curve with  $c_2$  depending only on 
$\lambda_0$, $c_0$ and $c_1$.

(1) We first prove $\ell(\gamma[a_0, x])\wedge \ell(\gamma[x, a_1])\le c\, d(x)$
 for all $x\in \gamma$ and some $c=c(c_1)$. 
Let $x_0\in \gamma $  be a point with maximal $d(x_0)$. By symmetry, it suffices to find an estimate  of the form 
$$\ell(\gamma[a_0, x])\le c\,d(x)$$
 for all  $x\in \gamma[a_0, x_0]$.   Let $q=q(c_1)\in (0, 1)$ be the 
number given by Lemma \ref{l3}.  If  $d(x)\le  q\, d(x_0)$, then   Lemma \ref{l3}  
implies $\ell(\gamma[a_0, x])\le M_2(c_1)\,d(x)$.  If  
$d(x)\ge  q\, d(x_0)$, we apply Lemma \ref{l2} with $r=d(x_0)$ 
 and obtain  $\ell(\gamma[a_0, x])\le \frac{M_1(c_1)}{q}\,d(x)$.\newline

(2) We next prove $\ell(\gamma)\le c_2\, d(a_0, a_1)$
  for some $c_2=c_2(\lambda_0, c_0, c_1)$.
 We may assume that $d(a_0)\le d(a_1)$.  
%Writing $t=d(a_0,a_1)$ we look for an estimate of the form $\ell(\gamma)\le c_2\,t$. 
Set $t=d(a_0,a_1)$ and  $r=d(a_0)$. 
 %By lemma ?? $\Omega\subset X$ is $(\lambda, c')$ quasiconvex.
We  consider two cases.

\noindent
{\bf{Case (a)}}.  $r\le 2c_0t/\lambda_0 $.   

We may assume  $\ell(\gamma)\ge 2t$.  Choose points 
 $b_0, b_1\in \gamma$ such that $\ell(\gamma[a_0,b_0])=\ell(\gamma[a_1,b_1])=t$.
 By  part (1) we have $t\le c\,d(b_i)$ for $i=0,1$.  We obtain:
$$ r_\Omega(b_0,b_1)\le \frac{d(b_0,a_0)+d(a_0, a_1)+d(a_1, b_1)}
{d(b_0)\wedge d(b_1)}\le \frac{3t}{t/{c}}=3c.$$

Since $\Omega$ is QH $c_1$-uniform, this implies that
$k(b_0, b_1)\le c_1 j_\Omega(b_0, b_1)\le c_1\log(1+3c)=c_3.$
 For each $x\in\gamma[b_0, b_1]$ we get  
$k(x, b_0)\le k(b_0, b_1)\le c_3$.  %By Lemma   \ref{bal}(2)  this implies that  
  %$d(x, b_0)\le d(b_0)(e^{c_3}-1)$. 
Since 
$d(b_0)\le d(a_0) +d(a_0, b_0)= r+t\le (1+2c_0/\lambda_0)t$,
  Lemma   \ref{bal}    yields
  $d(x)\le   %d(b_0)+d(x, b_0)\le 
d(b_0)e^{c_3}\le (1+2c_0/\lambda_0)te^{c_3}=c_4t$.  Integrating along  $\gamma[b_0, b_1]$ gives 
$k(b_0, b_1)\ge \ell(\gamma[b_0, b_1])/{c_4t}.$
 Since  $k(b_0, b_1)\le c_3$, we have
 $\ell(\gamma[b_0, b_1])\le c_3c_4 t$  and 
$\ell(\gamma)\le t+c_3c_4 t+t=(2+c_3c_4) t$.

\noindent
{\bf{Case (b)}}.  $r\ge 2c_0t/\lambda_0 $. 
%We  shall prove 

Since $\Omega\subset X$ is $(\lambda_0, c_0)$-quasiconvex  and 
$\frac{d(a_0, a_1)}{d(a_0)}=t/r\le \frac{\lambda_0 }{2c_0}$, 
  Lemma \ref{es} implies  $k(a_0, a_1)\le  2c_0\frac{t}{r}$.
%$d(a_1, a_0)=t\le \frac{\lambda_0 r}{2c_0}= 
%\frac{\lambda_0 d(a_0)}{2c_0}<\lambda_0 d(a_0)$,
 % there is a path $\omega$ from $a_0$ to $a_1$ with $\ell(\omega)\le c_0\,d(a_0, a_1)
%\le \lambda_0 r/2\le r/2$.  It follows that $d(x)\ge r/2$ for
 %all $x\in \omega$.
%Since $\delta(a_1)\ge \delta(a_0)=r$, we have 
 %Integrating   along $\omega$  gives 
%$k(a_0, a_1)\le \frac{2}{r}\ell(\omega)\le \frac{2}{r}c_0t$.
Set $b=\ell(\gamma)$ and let $\gamma':  [0,b]\ra \gamma$  be the arclength parametrization of $\gamma$  with  $\gamma'(0)=a_0$.  Since $d(\gamma'(s))\le d(a_0)+d(a_0, \gamma'(s))\le r+s$ for all $0\le s\le b$, we get 
$k(a_0, a_1)=\ell_k(\gamma)\ge \int^b_0\frac{ds}{r+s}=\log(1+b/r).$
 Now we have $\log(1+b/r)\le k(a_0, a_1)\le 2c_0\frac{t}{r}$. 
Setting $u=r/t\ge 2c_0/\lambda_0\ge 2$  we thus have 
$1+\frac{b}{tu}\le e^{2c_0/u}$.  It follows that $b/t\le c_6$, where
 $c_6=\max\{u(e^{2c_0/u}-1): u\ge 2\}$.  Notice that $c_6$  is finite
   and depends only on $c_0$.

\end{proof}

\begin{remark}
{Buckley and Herron (\cite{BH}) have recently  proved 
a result that  is  stronger   than 
Theorem \ref{cigart}.}
\end{remark}

\b{Th}\label{unifequ}
{Let $(X, d)$  be a proper metric space, and 
     $\Omega\subset X$ a $(\lambda_0, c_0)$-quasiconvex domain  for some
$0<\lambda_0\le 1$, $c_0\ge 1$.
%Suppose a proper domain $\Omega\subset X$ is $(\lambda_0, c_0)$ quasiconvex. 
  Then  the following   conditions  are quantitatively  
equivalent:\newline
(1) $\Omega$ is $c$-uniform;\newline
(2)  $\Omega$ is QH  $c$-uniform;\newline
(3) $k(x,y)\le c \,j_\Omega(x,y) +c'$   for
  all  $x, y\in \Omega$, where  $c$ and $c'$ are constants.}

\end{Th}

The phrase \lq\lq quantitatively  
equivalent" should be understood as follows. For example, 
 \lq \lq $(3)\Rightarrow (1)$ quantitatively"  means if 
$k(x,y)\le c \,j_\Omega(x,y) +c'$   for
  all  $x, y\in \Omega$,  then $\Omega$ is $c''$-uniform with 
$c''$ depending only on $c, c', \lambda_0$ and $c_0$.

\b{proof}
We show  that $(3)\Rightarrow(2)\Rightarrow(1)\Rightarrow(3)$.  The implication 
$(2)\Rightarrow(1)$ is simply Theorem \ref{cigart}.
  Assume that  (3) holds, let $a,b\in \Omega$, and set $r=r_\Omega(a,b)$.  If  
$r\le \frac{\lambda_0}{2c_0}$,  then $r\log 2\le \log(1+r)$,  and Lemma \ref{es}
implies that  $k(a,b)\le 2c_0r\le (2c_0/\log 2)j_\Omega(a,b)$.  If 
$r\ge \frac{\lambda_0}{2c_0}$, then $j_\Omega(a,b)\ge \log(1+\frac{\lambda_0}{2c_0})$. Hence 
$$\frac{k(a,b)}{j_\Omega(a,b)}\le c+\frac{c'}{\log(1+\frac{\lambda_0}{2c_0})},$$ 
   and we obtain (2).  It remains to prove $(1)\Rightarrow (3)$. 

 Assume that (1) is true, and let $a,b \in \Omega$.   
Let $\gamma\subset \Omega $  be a $c$-uniform arc  joining $a$ and $b$.
% and satisfying the uniformity conditions.   
Let $x_0\in \gamma$ be the point bisecting the length $\ell_0$ of $\gamma$. We may assume that $r_\Omega(a,b)\ge \frac{\lambda_0}{2c_0}$, since otherwise  Lemma 
\ref{es} gives 
$k(a,b)\le \lambda_0$ and we have $k(a,b)\le 0\cdot j_\Omega(a,b)+\lambda_0$. Setting
 $e=d(a)\wedge d(b)$ we have $\frac{e\lambda_0}{2c_0}\le d(a,b)\le \ell_0$. 
Choose points $a_1, b_1\in \gamma$ with
 $\ell(\gamma[a,a_1])=\ell(\gamma[b,b_1])=\frac{e\lambda_0}{4c_0}$.
 Notice that $r_\Omega(a, a_1), r_\Omega(b, b_1)\le \frac{\lambda_0}{2c_0}$.  Hence 
  Lemma \ref{es} yields $k(a_1, a)\le \lambda_0$ and  $k(b_1,b)\le \lambda_0$.  
Setting 
$\beta=\gamma[a_1, x_0]$ we obtain by the uniformity condition  
$$k(a_1, x_0)\le \ell_k(\beta)\le c\int_\beta\frac{ds(x)}{\ell(\gamma[a,x])}
=c\int^{\frac{\ell_0}{2}}_
\frac{e\lambda_0}{4c_0}\frac{ds}{s}=c\log\frac{2c_0\ell_0}{e\lambda_0}.$$
 The same  estimate also holds for $k(b_1, x_0)$. Since $\ell_o\le c \,d(a,b)$,
%Since a similar estimate holds for $k(b_1, x_0)$, 
  we have 
\begin{align*}
k(a,b) & \le k(a, a_1)+k(a_1,x_0)+k(x_0, b_1)+k(b_1, b)\\
  &\le 2\lambda_0+2c\log\frac{2c_0\ell_0}{e\lambda_0}\le 
2\lambda_0+2c\log\frac{2c_0cd(a,b)}{e\lambda_0}\\
  &=2\lambda_0+2c\log\frac{2c_0c}{\lambda_0}+2c\log\frac{d(a,b)}{e}\\
  &= c'+2c\log r_\Omega(a,b)
\le c'+2cj_\Omega(a,b),
\end{align*}
where $c'=2\lambda_0+2c\log\frac{2c_0c}{\lambda_0}$.
%.$$
 % Since $\log(Mt)\le M\log(1+t) $ for all $M\ge 1$, $t\ge 0$, this implies 
%$k(a, b)\le 4c^2j_\Omega(a,b)+1$.

\end{proof}

\subsection{Proof of the main results}\label{s4}

In this Section we prove the theorems stated in the Introduction.

%\b{Th}\label{t4.1}
%{Let $(X_i, d_i)$ ($i=1, 2$) be a proper metric space,
 %    $\Omega_i\subset X_i$ a domain with $\p\Omega_i\not=\emptyset$,  
  % and $\phi: \Omega_1\ra \Omega_2$
  %a   $\eta$-quasisymmetric map.  Suppose $(\Omega_1, d_1)$ is 
%$c_1$-uniform  and $(\Omega_2, d_2)$ is $(\lambda, c_2)$-quasiconvex.  Then 
%$(\Omega_2, d_2)$ is 
%$c$-uniform  with $c=c(\lambda, c_1, c_2, \eta)$. }

%\end{Th}
  
We first recall  two results of V\"ais\"al\"a.  % (Lemma 2.3 of \cite{V}).

\b{Le}\label{vle} (Lemma 2.3 of \cite{V})
{Suppose that $X$ is  $a$-quasiconvex,  $q>0$, $b\ge 0$, and that
$f:X\ra Y$ is a map with 
    $d(f(x), f(y))\le b$ whenever $d(x,y)\le q$.  Then 
$d(f(x), f(y))\le (ab/q)d(x,y)+b$ for all $x,y\in X$.}

\end{Le}

\b{Th}\label{6.12} (Theorem 6.12 of \cite{V})
{Suppose that $X$, $Y$ are metric spaces,  $A\subset X$,  $f:A\ra Y$
 is $\eta$-quasisymmetric,  and that $\overline{f(A)} $ is complete. Then
 $f$ extends to an $\eta$-quasisymmetric map $g: \ol A\ra Y$.}

\end{Th}

Let $L>0$  and  $A\ge 0$.  A map $f: X\ra Y$  between  two metric spaces is an 
  \e{$(L,A)$
quasi-isometry} if the following two conditions are satisfied:\newline
(1) $d(x_1,x_2)/L-A\le d(f(x_1), f(x_2))\le L\,d(x_1, x_2)+A$ holds for all 
$x_1, x_2\in X$;\newline
(2) For each $y\in Y$, there is some $x\in X$ with $d(f(x), y)\le A$.

\b{Le}\label{l4.2}
{For  $i=1,2$   
let $(X_i, d_i)$ %\e{($i=1,2$)} 
be a proper metric space,  
$\Omega_i\subset X_i$  a rectifiably connected domain
with $\p\Omega_i\not=\emptyset$,  and  
$g: \Omega_1\ra \Omega_2$
  an   $\eta$-quasisymmetric map.
Suppose $(\Omega_i, d_i)$ is $(\lambda_i, c_i)$-quasiconvex with $0<\lambda_i\le 1$ and
$c_i\ge 1$. 
 Let $k_i$ be the quasihyperbolic metric on
$\Omega_i\subset (X_i,  d_i)$.  Then 
the map $g: (\Omega_1, k_1)\ra (\Omega_2, k_2)$ 
 is  an  $(L,A)$ quasi-isometry with $L$ and $A$ depending only on 
 $\lambda_1$, $\lambda_2$, $c_1$,  $c_2$   and $\eta$.
}
\end{Le}

\b{proof}
By symmetry we only need to 
 show that there exist constants $L$ and $A$ depending only on
$\eta$, $\lambda_2$  and   $c_2$
such that $k_2(g(x), g(y))\le L \,k_1(x,y)+A$ for all $x,y\in \Omega_1$. 
 %The same argument yields an inequality
 %of the form $k_1(g^{-1}(x), g^{-1}(y))\le L'k_2(x,y)+A'$ for all $x,y\in \Omega_2$,
  %where $L'$ and $A'$ depend only on  $c''_1$ and $\eta''$.
 Since $(\Omega_1, k_1)$  is  a  %and $ (\Omega_2, k_2)$  are 
geodesic space,
by Lemma \ref{vle} it suffices to find a constant
 $q$ depending only on $\eta$, $\lambda_2$  and   $c_2$
 such that 
$k_2(g(x),g(y))\le \lambda_2$ whenever $k_1(x,y)\le q$.  
Let $q=\log[1+\eta^{-1}(\frac{\lambda_2}{2c_2})]$.
  Then  
  %We choose $q$ to be the unique number such that 
   $\eta(e^q-1)=\frac{\lambda_2}{2c_2}$.  %Since $\eta$ depends only on $\eta$,  
 Notice  that $q$ depends only on $\eta$, $\lambda_2$  and   $c_2$. We next
 show   that  $q$ has 
the required property.

Since $X_i$ are  proper  for $i=1,  2$   and 
$g: (\Omega_1, d_1)\ra (\Omega_2, d_2)$  is  $\eta$-quasisymmetric,
 Theorem \ref{6.12}  
 implies that 
$g$ extends to an $\eta$-quasisymmetric homeomorphism 
$(\overline\Omega_1, d_1)\ra (\overline\Omega_2, d_2)$, which is still denoted by $g$.
% between the closures of the domains.
Let  $x,y\in \Omega_1$ with $k_1(x,y)\le q$.
Then 
$$q\ge k_1(x,y)\ge \log\left(1+\frac{d_1(x,y)}{d_1(x)\wedge d_1(y)}\right)
\ge \log\left(1+\frac{ d_1(x,y)}{d_1(x)}\right),$$
 where $d_i(z)=d_i(z, \p\Omega_i)$ for  $z\in \Omega_i$. It follows that 
$\frac{d_1(x,y)}{d_1(x)}\le e^q-1$.
 Let $z\in\p\Omega_1$ with $d_2(g(x))=d_2(g(x), g(z))$.
Since $g$ is $\eta$-quasisymmetric, we have 
%\begin{align*}
\[
\frac{d_2(g(x), g(y))}{d_2(g(x))} 
=\frac{ d_2(g(x), g(y))}{ d_2(g(x), g(z))}
\le \eta\left(\frac{ d_1(x, y)}{ d_1(x, z)}\right) 
   \le \eta\left(\frac{ d_1(x, y)}{d_1(x)}\right)
\le \eta(e^q-1)=\frac{\lambda_2}{2c_2}.
\]
%\end{align*}
Since $(\Omega_2,  d_2)$ is $(\lambda_2, c_2)$-quasiconvex, 
Lemma \ref{es}  implies 
$ k_2(g(x), g(y))\le \lambda_2.$

%there is a path
 %$\gamma$ from $g(x)$ to $g(y)$ such that 
%$\ell(\gamma)\le c_2 \, d_2(g(x),g(y))$. 
 % It follows that $\ell(\gamma)\le c_2\, d_2(g(x),g(y))
%\le d_2(g(x))/2$ and hence $d_2(z)\ge d_2(g(x))/2$ for all
 %$z\in \gamma$.
%Now $$k_2(g(x), g(y))\le \int_\gamma \frac{1}{d_2(z)}|dz|\le 
%\frac{2}{d_2(g(x))}\ell(\gamma)\le 1.$$

\end{proof}

We also need  the following result  (Theorem 6.14 of \cite{V}).

\b{Th}\label{vai6.14}
{Suppose that $X$ is a connected metric space and that 
$f: X\ra Y$  is $\eta$-quasisymmetric. Then $f$ is $\eta_1$-quasisymmetric  for
     a  function  of the   form  
$\eta_1(t)=C(t^\alpha\vee t^{\frac{1}{\alpha}})$, where $C>0$  and $\alpha\in (0,1]$ depend only on $\eta$.}

\end{Th}

\b{Le}\label{l4.3}
{Let $(X_i, d_i)$ \e{($i=1, 2$)} be  proper metric spaces,
     $\Omega_i\subset X_i$  domains with $\p\Omega_i\not=\emptyset$,  
and $g: \Omega_1\ra \Omega_2$
  an   $\eta$-quasisymmetric map.
Then there are constants $a>0$ and $b>0$ depending only on $\eta$  such that  the 
 following hold:
\[
 j_{\Omega_2}(g(x_1), g(y_1))\le a\, j_{\Omega_1}(x_1, y_1) + b\;\;\;
 {\text{for all}} \;\;\; x_1, y_1\in \Omega_1
\] and
 \[
 j_{\Omega_1}(g^{-1}(x_2), g^{-1}(y_2))\le a \,j_{\Omega_2}(x_2, y_2) + b\;\;\;
 {\text{for all}}   \;\;\; x_2, y_2\in \Omega_2.
\]
}
\end{Le}

\b{proof}
By symmetry it suffices to prove that there 
are constants $a$ and $b$ depending only on $\eta$  such that
$ j_{\Omega_1}(g^{-1}(x_2), g^{-1}(y_2))\le a\, j_{\Omega_2}(x_2, y_2) + b$
 for all $x_2, y_2\in \Omega_2$.   To do so,  we shall 
%it suffices  to 
      find constants 
 $c>1$ and $d\ge 1$ depending only on $\eta$
such that 
  $1+r_{\Omega_1}(g^{-1}(x_2), g^{-1}(y_2))\le 
 c\,(1+r_{\Omega_2}(x_2,  y_2))^d$   for all $x_2, y_2\in \Omega_2$.

%$1+r_1\le c(1+r_2)^d$ for all $x_2, y_2\in \Omega_2$, where
%$$r_1=\frac{d_1(g^{-1}(x_2), g^{-1}(y_2))}
%{d_1(g^{-1}(x_2))\wedge d_1(g^{-1}(y_2))}\;\;\; {\text{ and}}\;\;\;
%r_2=\frac{d_2(x_2, y_2)}{d_2(x_2)\wedge d_2(y_2)}.$$

%Since $X_i$ is complete, 
%Theorem 6.12
%of \cite{V} implies that 
  By  Theorem \ref{6.12},   
 %As in the proof of Lemma \ref{l4.2}, 
$g$ extends to a $\eta$-quasisymmetric homeomorphism 
$\bar{g}: (\overline\Omega_1, d_1)\ra (\overline\Omega_2, d_2)$.
 Then (see Theorem 6.3 of \cite{V})
  $\bar{g}^{-1}: (\overline\Omega_2, d_2)\ra (\overline\Omega_1, d_1)$
is $\eta'$-quasisymmetric with $\eta'(t)=\eta^{-1}(t^{-1})^{-1}$.
By Theorem \ref{vai6.14}, %we may assume that 
   the map $\bar g^{-1}$ is $\eta_1$-quasisymmetric with
$\eta_1(t)=C(t^\alpha\vee t^{\frac{1}{\alpha}})$, where $C>0$  and $\alpha\in (0,1]$ depend only on $\eta$.
Set $c=1+C$ and 
$d=\frac{1}{\alpha}$.

Fix $x_2, y_2\in \Omega_2$ and 
   set $r_1=r_{\Omega_1}(g^{-1}(x_2), g^{-1}(y_2))$,  
$r_2=r_{\Omega_2}(x_2,  y_2)$.    
  %define $r_1, r_2$ as above. 
We may assume $d_1(g^{-1}(x_2))\le d_1(g^{-1}(y_2))$. Pick $w\in \p \Omega_1$
  with $d_1(g^{-1}(x_2))=d_1(g^{-1}(x_2), w)$.  Since $\bar g(w)\in \p \Omega_2$, we have 
$d_2( x_2, \bar g(w))\ge d_2(x_2)\ge d_2(x_2)\wedge d_2(y_2)$. 
Now 
\begin{align*}
 r_1 & =\frac{d_1(g^{-1}(x_2), g^{-1}(y_2))}
{d_1(g^{-1}(x_2))\wedge d_1(g^{-1}(y_2))}
%\frac{d_1(x_1, y_1)}{d_1(x_1)\wedge d_1(y_1)}
=\frac{d_1(g^{-1}(x_2), g^{-1}(y_2))}{d_1(g^{-1}(x_2), w)}\\
   & \le \eta_1\left(\frac{d_2(x_2, y_2)}
{d_2(x_2, \bar g(w))}\right)
\le \eta_1\left(\frac{d_2(x_2, y_2)}{d_2(x_2)\wedge d_2(y_2)}\right)
=\eta_1(r_2),
\end{align*}
  that is,  $r_1\le \eta_1(r_2)$. 

If $r_2\le 1$, then $r_1\le \eta_1(r_2)\le \eta_1(1)=C$, hence 
$1+r_1\le 1+C=c\le c(1+r_2)^d$.   If $r_2\ge 1$,  
  then $r_1\le \eta_1(r_2)=Cr_2^d$. It follows that
$1+r_1\le 1+Cr_2^d\le c(1+r_2^d)
\le c(1+r_2)^d$.

\end{proof}

\noindent
{\bf{Proof of Theorem \ref{t4.1}.}}
The assumptions imply that for  $i=1,  2$,      
$(\Omega_i, d_i)$ is rectifiably connected,
 $(\Omega_1, d_1)$  is $(1/2, c_1)$-quasiconvex
 and $(\Omega_2, d_2)$  is 
$(\lambda, c_2)$-quasiconvex. 
%Notice the asumptions also imply that 
Fix $x_2, y_2\in \Omega_2$. By Lemma \ref{l4.2}, $k_2(x_2, y_2)\le
L \,k_1(g^{-1}(x_2), g^{-1}(y_2))+A$, where 
$L=L(\eta, \lambda,  c_1, c_2)$ and $A=A(\eta, \lambda, c_1, c_2)$. 
Since $(\Omega_1, d_1)$ is $c_1$-uniform, Theorem  \ref{unifequ}  implies 
$$k_1(g^{-1}(x_2), g^{-1}(y_2))\le c'\, j_{\Omega_1}(g^{-1}(x_2), g^{-1}(y_2)),$$ where 
$c'=c'(c_1)$. 
  On the other hand,  by Lemma \ref{l4.3}, 
\[j_{\Omega_1}(g^{-1}(x_2), g^{-1}(y_2))\le a\, j_{\Omega_2}(x_2, y_2)+b,  \]
 where
$a=a(\eta)$, $b=b(\eta)$. 
Combining the above inequalities we have
$k_2(x_2, y_2)\le a'\, j_{\Omega_2}(x_2, y_2)+b'$ for all $x_2, y_2\in \Omega_2$,
 where $a', b'$ depend only on $\eta$, $\lambda$, $c_1$ and $c_2$.
Since $(\Omega_2, d_2)$ is $(\lambda, c_2)$-quasiconvex, Theorem \ref{unifequ}
 implies that
$(\Omega_2, d_2)$  is $c$-uniform with $c=c(\eta, \lambda, c_1, c_2)$.

\qed

In the proof of the following lemma, we shall implicitly 
use the inequality
$\frac{1}{4}s_p(x,y)\le \hat{d_p}(x,y)\le s_p(x,y)$  %implicitly
(see Section \ref{s2}).

\b{Le}\label{ll}
{Let $(X, d)$ be an unbounded  proper metric space,
 $\Omega\subset X$ an unbounded  domain,   $p\in \p\Omega$,
      and   $0<\lambda\le 1$, $c\ge 1$.
  % be constants. 
If $\Omega\subset (X, d)$ is $(\lambda, c)$-quasiconvex, then
 $\Omega\subset (X, \hat{d}_p)$  is $(\lambda', c')$-quasiconvex
  with $\lambda'=\lambda'(\lambda, c)$ and $c'=c'(\lambda, c)$.
}

\end{Le}

\b{proof}
Set $\lambda'=\frac{\lambda}{10000c^2}$ and $c'=64c$.
  Fix  $x\in \Omega$. Let $\hat{d}_p(x)=\hat{d}_p(x,\hat{\p}\Omega)$
  denote  the $\hat{d}_p$-distance from $x$ to $\hat{\p}\Omega$, 
where $\hat{\p}\Omega=\p\Omega\cup \{\i\}$ 
is the boundary of $\Omega$ 
in $(X\cup \{\i\}, \hat{d}_p)$.   %  and $\hat{d}_p(x)=\hat{d}_p(x,\hat{\p}\Omega)$
   %is the $\hat{d}_p$-distance from $x$ to $\hat{\p}\Omega$. 
For $r>0$, let $\hat{B}_p(x,r):=\{y\in S_p(X):  \hat{d}_p(x,y)< r\}$.

%Notice that $\hat{\p}\Omega=\p\Omega$ if $(\Omega, d)$ is bounded and 
%$\hat{\p}\Omega=\p\Omega\cup \{\i\}$  if $(\Omega, d)$ is unbounded.
%For $r>0$, let $\hat{B}_p(x,r):=\{y\in S_p(X):  \hat{d}_p(x,y)< r\}$.

\noindent
{\bf{Claim:}} $\hat{B}_p(x, \lambda'\hat{d}_p(x))\subset B(x, \frac{\lambda}{10c}d(x))$.

We first assume the claim and complete  the proof of the lemma.
Let $y_1, y_2\in \hat{B}_p(x, \lambda'\hat{d}_p(x))$. 
Since $\Omega\subset (X, d)$ is $(\lambda, c)$-quasiconvex,
 the claim implies that there is a path $\gamma$ from $y_1$ to $y_2$ such that
$\ell(\gamma)\le c \, d(y_1, y_2)$. The claim further implies 
$d(y, x)\le 3d(x)/10$ for all $y\in \gamma$. 
Since $p\in \p\Omega$, we have  $d(x,p)\ge d(x)$. It follows that 
 for any $y\in \gamma$   we have 
$d(x,p)/2\le d(y,p)\le 2 d(x,p)$.    Hence
$(1+d(x,p))/2\le 1+d(y,p)\le 2(1+d(x,p))$. 
Because 
$$\frac{d(z_1, z_2)}{4(1+d(z_1,p))(1+d(z_2,p))}\le \hat{d}_p(z_1, z_2)
\le \frac{d(z_1, z_2)}{(1+d(z_1,p))(1+d(z_2,p))}$$
 for all $z_1, z_2\in X$,
we conclude that 
$$\hat{\ell}_p(\gamma)\le \frac{4\ell(\gamma)}{(1+d(x,p))^2}\;\;\;{\text{and}}
 \;\;\; \hat{d}_p(y_1,y_2)\ge \frac{d(y_1, y_2)}{16(1+d(x,p))^2},$$
  where $\hat{\ell}_p(\gamma)$ denotes the $\hat{d}_p$-length of $\gamma$.
 Together  with $\ell(\gamma)\le c \; d(y_1, y_2)$, the above two inequalities imply 
$$\hat{\ell}_p(\gamma)\le 64c \hat{d}_p(y_1,y_2)=c'\hat{d}_p(y_1,y_2).$$
We have shown that $\Omega$ is $(\lambda', c')$-quasiconvex
in $(X, \hat{d}_p)$. 

Next we prove the claim. Let $y\in \hat{B}_p(x, \lambda'\hat{d}_p(x))$.
We need to prove $d(x,y)< \frac{\lambda}{10c}d(x)$.
  There is some $w\in \p\Omega$ with $d(x)=d(x, w)$
 and some $z\in \hat{\p}\Omega$ with $\hat{d}_p(x)=\hat{d}_p(x,z)$.
  % where $\hat{\p}\Omega=\p\Omega\cup \{\i\}$ 
  %is the boundary of $\Omega$ 
 %in $(X\cup \{\i\}, \hat{d}_p)$  and $\hat{d}_p(x)=\hat{d}_p(x,\hat{\p}\Omega)$
 %is the $\hat{d}_p$-distance from $x$ to $\hat{\p}\Omega$. 
We consider two cases depending on whether $z=\i$ or not.

\noindent
{\bf{Case (1).}} $z=\i$. Then $\hat{d}_p(x)=\hat{d}_p(x,\i)\le \frac{1}{1+d(x,p)}$.
The fact that 
$y\in \hat{B}_p(x, \lambda'\hat{d}_p(x))$ now 
implies $d(x,y)\le 4\lambda'(1+d(y,p))$. 
Since $w\in \p\Omega$, we have $\hat{d}_p(x,w)\ge \hat{d}_p(x,\i)$,
 which implies that   $d(x,w)\ge (1+d(w,p))/4$. 
If  $d(y,p)\le 1$,  
 then 
\begin{align*}
d(x,y) & \le 4\lambda'(1+d(y,p))\le 8\lambda'=\frac{8\lambda}{10000c^2}
< \frac{\lambda}{40c}\\
  & \le\frac{\lambda}{10c} \frac{1+d(w,p)}{4}\le 
 \frac{\lambda}{10c}d(x,w)=\frac{\lambda}{10c}d(x),
\end{align*}
 and we are done.  Now assume $d(y,p)\ge 1$.
  In this case 
 we  shall prove $d(x,y)< \frac{\lambda}{10c}d(x)$ by contradiction. 
So we suppose $d(x,y)\ge \frac{\lambda}{10c}d(x)$.
%First assume $d(y,p)\ge 1$.  
 We have  $d(x,y)\le 4\lambda'(1+d(y,p))\le 8\lambda'd(y,p)$.
Hence $d(y,p)/2\le d(x,p)\le 2\,d(y,p)$  and 
$$d(x,w)=d(x)\le d(x,y)\frac{10c}{\lambda}\le 
8\lambda'd(y,p)\frac{10c}{\lambda}\le \frac{d(y,p)}{100}\le \frac{d(x,p)}{50}.$$
It follows that $d(w,p)\ge d(x,p)/2$,   and  therefore 
$$d(x,w)\ge (1+d(w,p))/4\ge d(w,p)/4\ge d(x,p)/8,$$ contradicting
 $d(x,w)\le d(x,p)/50$.

%Next we assume $d(y,p)\le 1$.
 %Then 
%\begin{align*}
%d(x,y) & \le 4\lambda'(1+d(y,p))\le 8\lambda'=\frac{8\lambda}{10000c^2}
 %< \frac{\lambda}{40c}\\
  %& \le\frac{\lambda}{10c} \frac{1+d(w,p)}{4}\le 
 %\frac{\lambda}{10c}d(x,w)=\frac{\lambda}{10c}d(x),
%\end{align*}
 %contradiction again. 

\noindent
{\bf{Case (2).}}  $z\in \p\Omega$.  
The inequalities  $\hat{d}_p(x,y)\le \lambda' \hat{d}_p(x,z)$,
$\hat{d}_p(x,z)\le \hat{d}_p(x,w)$ and $\hat{d}_p(x,z)\le \hat{d}_p(x,\i)$ 
yield the inequalities 
 $$\frac{d(x,y)}{4(1+d(y,p))}\le \lambda'\frac{d(x,z)}{1+d(z,p)},\;
 \frac{d(x,z)}{(1+d(z,p))}\le \frac{4 d(x,w)}{1+d(w,p)} \;\;{\text{and}}\;\;
 \frac{d(x,z)}{(1+d(z,p))}\le 4$$ respectively.
  If  $d(y,p)\le 1$,  then 
\begin{align*}
d(x,y)  &  \le 4(1+d(y,p))\cdot 
  \lambda'\frac{d(x,z)}{1+d(z,p)}
\le 8\lambda'\frac{4d(x,w)}{1+d(w,p)}\\
   & \le 32\lambda'd(x,w)
< 
\frac{\lambda}{10c} d(x,w)=\frac{\lambda}{10c} d(x).
\end{align*}

 On the other hand,  if   $d(y,p)\ge 1$, then 
   $$d(x,y)\le 4(1+d(y,p)) \cdot \lambda'\frac{d(x,z)}{1+d(z,p)}
\le 8d(y,p)\cdot\lambda'\cdot 4=32\lambda'd(y,p)\le d(y,p)/10.$$
%It follows that 
 If $d(x,w)\ge \frac{d(y,p)}{10c}$, then 
   $d(x,y)\le 32\lambda'd(y,p)\le 32\lambda' 10cd(x,w)<\frac{ \lambda}{10c} d(x,w)$, 
 and we are done. 
 If $d(x,w)\le \frac{d(y,p)}{10c}$, then 
$d(w,y)\le d(w,x)+d(x,y)\le d(y,p)/5$ and   hence
$d(w,p)\ge d(y,p)/2$.  It follows that
\begin{align*}
d(x,y)  &  \le 4(1+d(y,p))\cdot\lambda'\frac{d(x,z)}{1+d(z,p)}
\le 8d(y,p)\cdot\lambda' \cdot \frac{4 d(x,w)}{1+d(w,p)}\\
  &  \le 32\lambda'd(y,p)\frac{d(x,w)}{d(w,p)}
\le 32\lambda' \cdot 2d(x,w)=64\lambda'd(x,w)\\
  & < 
\frac{\lambda}{10c} d(x,w)=\frac{\lambda}{10c} d(x).
\end{align*}

\end{proof}

\noindent
{\bf{Proof of Theorem \ref{main}.}}
 If $(\Omega_i, d_i)$ ($i=1, 2$) is bounded, let $d'_i=d_i$; and if 
$(\Omega_i, d_i)$  is unbounded,
 pick $p_i\in \p\Omega_i$  and 
 set  $d'_i:=\hat{d_i}_{p_i}$.  Recall that the identity map
$(X_i, d_i)\ra (X_i, d'_i)$ is quasim\"obius. 
  As  $(\Omega_i, {d'_i})$
 is bounded,  a
%Since a quasim\"obius map between bounded metric spaces is quasisymmetric, 
 %Theorem \ref{main} follows from Theorem \ref{t4.1} when both
 %$(\Omega_1, d_1)$ and  $(\Omega_2, d_2) $ are bounded.
 %In general, we consider the metric  $\hat{d_i}_{p_i}$ ($p_i\in X_i$)
  % on $X_i$ ($i=1, 2$),
 %then a 
quasim\"obius map 
$g: (\Omega_1, d_1) \ra (\Omega_2, d_2)$  becomes a quasisymmetric 
map $g: (\Omega_1, {d'_1}) \ra (\Omega_2, {d'_2})$.
% as
 %$(X_i, \hat{d_i}_{p_i})$
 %is bounded.  
Since $(\Omega_1, d_1)$ is uniform, Theorem \ref{bk2}
 implies that  $(\Omega_1, {d'_1})$ is uniform. On the other hand,
$(\Omega_2, d_2)$ is $(\lambda, c_2)$-quasiconvex.  By 
Lemma \ref{ll},  $(\Omega_2, {d'_2})$
is $(\lambda', c')$-quasiconvex for some $0<\lambda'\le 1$ and $c'\ge 1$.
Now  it follows from 
Theorem \ref{t4.1}
  that $(\Omega_2, {d'_2})$ is uniform.  Hence 
$(\Omega_2, d_2)$ is also uniform by Theorem \ref{bk2}. 

\qed

\noindent
{\bf{Proof of Theorem \ref{ac} (2).}}
Let  $i=1$  or  $ 2$.  If $\Omega_i$ is bounded, 
set $X'_i=X_i$ and  $d'_i=d_i$; %and $f_i: (X_i,  d_i)\ra (X'_i,  d'_i)$
  %the identity map;
  if $\Omega_i$ is unbounded,  then fix any base point $p_i\in \p\Omega_i$ and set
$X'_i=S_{p_i}(X_i)$  and  $d'_i=\hat{d_i}_{p_i}$.  % here $d=d_i$. 
Denote by $\p\Omega'_i$ the boundary of $\Omega_i$ in $(X'_i, d'_i)$
 and  $\overline\Omega'_i$ the  closure  of $\Omega_i$ in $(X'_i, d'_i)$.
 Let $f_i: (\Omega_i,  d_i)\ra (\Omega_i,  d'_i)$
 be the identity map and set 
$$h':=f_2\circ h\circ f_1^{-1}: (\Omega_1,  d'_1)\ra (\Omega_2,  d'_2).$$  %Also let $\Omega'_i=(\Omega_i, d'_i)$ and 
%Denote by $\p\Omega'_i$ the boundary of $\Omega_i$ in $(X'_i, d'_i)$
 %and  $\overline\Omega'_i$ the  closure  of $\Omega_i$ in $(X'_i, d'_i)$.
Let $\eta_0(t)=16t$.
Then $h'$ is a $\eta':=\eta_0\circ \eta\circ \eta_0$-quasim\"obius 
   homeomorphism between bounded metric spaces,
 and hence is a quasisymmetric map.
By Theorem \ref{6.12},
  the map $h'$ extends continously to a  
%$\eta'$-quasim\"obius 
homeomorphism
   $(\overline\Omega'_1, d'_1)\ra (\overline\Omega'_2, d'_2)$, which is still denoted by $h'$. 
In particular, 
there exist $a_1\in \p\Omega'_1$,  $a_2\in \p\Omega'_2$
  such that for any $\{x_i\}\subset  \Omega_1$ with $x_i\ra a_1$ we have $h'(x_i)\ra a_2$.
If $\p\Omega'_1$ is a single point, then $\p\Omega'_2$ and $\p\Omega_2$ are also single 
points.  Since $(X_2, d_2)$ is $c_2$-quasiconvex and $c_2$-annular convex, 
the fact that $\p\Omega_2$  is a single point implies that 
%it is not hard to check that 
$(\Omega_2, d_2)$ is $6c^2_2$-uniform
%In this case  one checks that $(\Omega_2, d_2)$ is $6c^2_2$-uniform 
(see Lemma 9.4 of \cite{HSX}).
%By Lemma \ref{onepo} $(\Omega_2, d_2)$ is $6c^2_2$-uniform.
 From now on, we assume  $\p\Omega'_1$ contains at least two points.

Now we fix $a_1\in \p\Omega'_1, a_2\in \p\Omega'_2$ 
such that for any $\{x_i\}\subset  \Omega_1$ with $x_i\ra a_1$ we have $h'(x_i)\ra a_2$.
Let $X''_i=I_{a_i}(X'_i)=X'_i\backslash \{a_i\}$,  $d''_i=(d'_i)_{a_i}$.
Denote by $\p\Omega''_i$ the boundary of $\Omega_i$ in $(X''_i, d''_i)$
 and  $\overline\Omega''_i$ the  closure  of $\Omega_i$ in $(X''_i, d''_i)$.
Note   that    $\p\Omega''_i=\p\Omega'_i\backslash \{a_i\}$
 and $\overline\Omega''_i=\overline\Omega'_i\backslash \{a_i\}$ as sets. 
 Let $g_i: (X'_i\backslash \{a_i\},  d'_i)\ra (X'_i\backslash \{a_i\},  d''_i)$
 be the identity map and set 
$$h'':=g_2\circ h'\circ g_1^{-1}: 
(\overline\Omega''_1,  d''_1)\ra (\overline\Omega''_2,  d''_2).$$  %Also let $\Omega'_i=(\Omega_i, d'_i)$ and 
Since $g_i$ is $\eta_0$-quasim\"obius,   
$h''$  is $\eta''$-quasim\"obius, 
where $\eta'':=\eta_0\circ \eta'\circ \eta_0$.
 The choice of $a_1$ and $a_2$ implies that for any 
$\{x_i\}\subset  \Omega_1$ with $d''_1(x_i, x_1)\ra \i$
we have $d''_2(h''(x_i), h''(x_1))   \ra \i$.
% Lemma \ref{l1}  implies that 
%for any $\{x_i\}\subset  \Omega_1$ with $x_i\ra \i$ we have $h'(x_i)\ra \i$.
It follows that $h''$  is an  $\eta''$-quasisymmetric homeomorphism.

Since $(\Omega_1, d_1)$ is $c_1$-uniform,
Theorem \ref{bk2}  (1) implies that $(\Omega_1, d'_1)$ is $c'_1$-uniform
 with $c'_1=c'_1(c_1)$.  %Theorem \ref{spe}  (3) and the fact that 
%$(X_1,d_1)$ is $c_1$-quasiconvex and $c_1$-annulus quasiconvex  imply that
%$(X'_1,d'_1)$ is $c'_1$-quasiconvex and $c'_1$-annulus quasiconvex
 %with $c'_1=c'_1(c_1)$.   
Since $\p\Omega'_1$  contains at least two points and 
$a_1\in \p\Omega'_1$, it follows from 
Theorem \ref{bk3} (1)  that 
$(\Omega_1, d''_1)$ is $c''_1$-uniform 
with  $c''_1=c''_1(c'_1)=c''_1(c_1)$.  
%The fact that $(\Omega_1, d''_1)$ is $c''_1$-uniform  implies that 
%$(\overline\Omega_1'', d''_1)$ is $c''_1$-quasiconvex.  
On the other hand, 
since  $(X_2,d_2)$ is $c_2$-quasiconvex and $c_2$-annular convex,
  %Theorem \ref{spe} (3)   and  
 %Theorem \ref{inve} (1) 
it follows from Section \ref{s2}  that 
$(X''_2,d''_2)$ is $c''_2$-quasiconvex % and $c''_2$-annulus quasiconvex
    with
 $c''_2=c''_2(c_2)$. 
  Therefore,    $\Omega_2\subset  (X''_2,d''_2)$ is 
$(1/(3c''_2), c''_2)$-quasiconvex.

Now Theorem  \ref{t4.1}  applied to $h''$
 %and  \ref{bdcase}
 implies  that $(\Omega_2, d''_2)$ is $c'$-uniform
 with $c'=c'(1/(3c''_2),  c''_2, c''_1,  \eta'')=c'(c_1, c_2, \eta)$. 
  Now the result follows from   Theorem \ref{bk3}  (2) and Theorem \ref{bk2} (2).

\qed

\noindent
{\bf{Proof of Theorem \ref{ac} (1).}}
   The proof is  similar to  that  of  Theorem \ref{ac}  (2),   
  %We repeat the proof of Theorem \ref{ac}, 
and we only indicate what should be modified.
By the assumption of Theorem \ref{ac} (1) $(\Omega_2, d_2)$ is unbounded,
hence  $\i\in \p\Omega'_2$.  We choose $a_2=\i$ and $a_1=(h')^{-1}(\i)$.
 % and so $X'_2=X_2$  and $d'_2=d_2$.  
The proof of Theorem \ref{ac}  (2)  shows that 
$(\Omega_1, d''_1)$ is $c''_1$-uniform 
with  $c''_1=c''_1(c_1)$  
  and  $h''$  is an  $\eta''$-quasisymmetric homeomorphism
  with $\eta''=\eta''(\eta)$.

Since $a_2=\i$, Lemma \ref{rebm1}
  implies that  the identity map $(X_2, d_2)\ra  (X_2, d''_2)$ 
 is $16$-bilipschitz. 
Now the fact that $(\Omega_2, d_2)$ is $(\lambda, c_2)$-quasiconvex implies that
 $(\Omega_2, d''_2)$ is $(\lambda'', c''_2)$-quasiconvex
 with $\lambda''=\lambda/256$ and $c''_2=256 c_2$.

Now Theorem  \ref{t4.1}  applied to $h''$
 %and  \ref{bdcase}
 implies  that $(\Omega_2, d''_2)$ is $c'$-uniform
 with $c'=c'(\lambda'',  c''_2, c''_1,  \eta'')=c'(\lambda, c_1, c_2, \eta)$. 
  Now the result follows from   Lemma \ref{rebm1}.

%Theorem \ref{bk3}  (2) and Theorem \ref{bk2} (2). 

\qed

\subsection{Example and open questions}\label{s5}

In this Section we give an example that shows Theorem \ref{main}  can not be made 
quantitative,  and  present   two related questions.

%noindent
%\bf{Example}}.

Let $\R^n$ be the $n$-dimensional Euclidean space  and 
$\tau: \R^n\backslash\{p\}\ra \R^n\backslash\{p\}$,  $\tau(x)={x}/{|x|^2}$ 
  %$\tau(x)=\frac{x}{|x|^2}$ 
the inversion 
about the unit sphere centered at  the origin $p$. 
Let $d$ denote the Euclidean metric. We can define a new metric 
$d'$  on
 $\R^n\backslash\{p\}$  by pulling back the Euclidean metric via $\tau$:
$d'(x,y)=d(\tau(x), \tau(y))$.  One checks that
%$$\leqno{(*)}\;\;\;\;
$$ d'(x,y)=\frac{d(x,y)}{d(x,p)d(y,p)}.  $$
 Let $A\subset \R^n$  be  a subset containing $p$,
  and consider the metric spaces  $(A, d)$  and $(A\backslash\{p\}, d_p)$.   
  Notice that  $f_p(x,y)=d'(x,y)$ for all $x,y\in A\backslash\{p\}$.
%n other words,  $f_o(x,y)=d'(x,y)$ for all $x,y\in \R^n\backslash\{o\}$.
 Since $d'$ is a metric on $\R^n\backslash\{p\}$,  %or any subset 
%A\subset \R^n\backslash\{o\}$, 
$f_p$ is a
metric on $A\backslash\{p\}$. 
 Now the definition of $d_p$ and the triangle
  inequality show that  $d_p=f_p$  on  $A\backslash\{p\}$.   It follows that 
for any $x,y\in A\backslash\{p\}$, we have 
$d_p(x,y)=f_p(x,y)=d'(x,y)=d(\tau(x), \tau(y))$;
  that is,   
$\tau:  (A\backslash\{p\}, d_p)\ra  (\tau(A\backslash\{p\}), d)$ is  
an isometry.

Now consider $\R^2$.  We identity $\R^2$ with $\C$ and 
use complex number notations. For $0<u<\pi/2$, let 
$X=\{\frac{1}{2}(i+e^{i\theta}):  -\pi/2\le \theta \le 3\pi/2-u\}$
  and $\Omega=X\backslash\{p, q\}$,  where $q=\frac{1}{2}(i+e^{i(3\pi/2-u)})$.
  One checks that $\Omega\subset (X, d)$  is  a  $(1/2, \pi)$-quasiconvex
 domain.   By Section \ref{s2},  the identity map
$(X\backslash\{p\}, d)\ra (X\backslash\{p\}, d_p)$  is $\eta$-quasim\"obius
  with $\eta(t)=16t$.
By the preceding paragraph,
 $(X\backslash\{p\}, d_p)$   is isometric to  $(\tau(X\backslash\{p\}), d)$.
Set $u'=\frac{\cos(3\pi/2-u)}{1+\sin(3\pi/2-u)}$. We notice that
$\tau(X\backslash\{p\})=\{x+i:  u'\le x<\i\}$  is a ray
  and $\tau(\Omega)=\{x+i:  u'<
 x<\i\}$.  It is now 
clear  that $\Omega\subset (X\backslash\{p\}, d_p)$   is $1$-uniform. 
  On the other hand,  by considering  two points in $\Omega$ close to 
 $p$ and $q$, we see that 
the uniformity constant
  of $\Omega\subset (X, d)$  is in the order of $1/u$, which tends to infinity as $u\ra 0$.
  This example shows that 
Theorem \ref{main}  can not be made 
quantitative.

In view of the above example and the Theorems in this paper, 
it is natural to ask the following question:

\b{question}\label{Q1}
{Let $(X_i, d_i)$ \e{(}$i=1, 2$\e{)} be a proper metric space,  and 
     $\Omega_i\subset X_i$ a domain with $\p\Omega_i\not=\emptyset$.  
  Suppose $(\Omega_1, d_1)$ is $c_1$-uniform,    
$(\Omega_2, d_2)$ is   bounded  and  $c_2$-quasiconvex,  
% for some $0<\lambda\le 1$ and $c\ge 1$,   
and   
    there is an  $\eta$-quasim\"obius homeomorphism
$ (\Omega_1, d_1)\ra (\Omega_2, d_2)$.
By Theorem \ref{main},   $(\Omega_2, d_2)$ is $c$-uniform
  for  some constant $c$. 
 Is it possible   to obtain an upper  bound for $c$
in terms only  of 
   $c_1$, $c_2$ and $\eta$?
}

\end{question}

A special case of Question \ref{Q1} is the following:

\b{question}\label{Q2}
{Let $(X,d)$ be a proper metric space,
$\Omega\subset X$ a bounded domain and $p\in \p\Omega$.
Assume $\p\Omega$ contains at least two points, 
$(\Omega, d_p)$ is $c_1$-uniform
   and  $(\Omega, d)$ is $c_2$-quasiconvex.
%  and $(\Omega, d_p)$ is $c'$-uniform. 
 By Theorem \ref{bk3},   $(\Omega, d)$ is $c$-uniform
   for  some  constant  $c$. 
 Is it possible   to obtain an upper  bound for $c$
in terms only  of 
   $c_1$  and $c_2$ ?

%true  that $c''$ can be bounded from above in terms of 
 %  $c$, $c'$ and $\lambda$?

}

\end{question}

 \addcontentsline{toc}{subsection}{References}
%\noindent 
%Xiangdong Xie\\
%Department of Mathematics\\
%Washington University \\
%St.Louis, MO 63130\\
%xxie@math.wustl.edu\\

\end{document}